\input amstex\documentstyle{amsppt}  
\pagewidth{12.5cm}\pageheight{19cm}\magnification\magstep1
\topmatter
\title Total positivity in Springer fibres\endtitle
\author G. Lusztig\endauthor
\address{Department of Mathematics, M.I.T., Cambridge, MA 02139}\endaddress
\thanks{Supported by NSF grant DMS-1855773.}\endthanks
\endtopmatter   
\document

\define\pos{\text{\rm pos}}

\define\dw{\dot w}

\define\ds{\dot s}

\define\mpb{\medpagebreak}

\define\frl{\forall}

\define\si{\sim}

\define\sqc{\sqcup}

\define\qua{\quad}

\define\op{\oplus}
   
\define\part{\partial}
\define\emp{\emptyset}

\define\n{\notin}
\define\iy{\infty}
\define\m{\mapsto}
\define\do{\dots}

\define\lra{\leftrightarrow}

\define\sm{\smallmatrix}
\define\esm{\endsmallmatrix}
\define\sub{\subset}    

\define\T{\times}
\define\ti{\tilde}
\define\nl{\newline}
\redefine\i{^{-1}}
\define\fra{\frac}
\define\un{\underline}

\define\Ad{\text{\rm Ad}}
\define\Hom{\text{\rm Hom}}

\define\supp{\text{\rm supp}}

\define\a{\alpha}
\redefine\b{\beta}
\redefine\c{\chi}
\define\g{\gamma}
\redefine\d{\delta}
\define\e{\epsilon}
\define\et{\eta}

\define\p{\pi}
\define\ph{\phi}

\define\s{\sigma}
\redefine\t{\tau}

\define\z{\zeta}
\define\x{\xi}

\define\Ph{\Phi}

\define\ii{\bold i}

\define\CC{\bold C}

\define\NN{\bold N}

\define\RR{\bold R}

\define\cb{\Cal B}

\define\ci{\Cal I}

\define\cl{\Cal L}

\define\cp{\Cal P}

\define\cu{\Cal U}

\define\cz{\Cal Z}
\define\cx{\Cal X}
\define\cy{\Cal Y}

\define\fb{\frak b}

\define\fg{\frak g}

\define\fp{\frak p}

\define\fL{\frak L}

\define\fU{\frak U}

\define\fZ{\frak Z}
\define\fX{\frak X}

\head Contents\endhead
1. The totally positive part of $\cb_u$.

2. Examples.

3. Partial flag manifolds.

4. A conjecture and its consequences.

5. The map $g\m P_g$ from $G_{\ge0}$ to $\cp_{\ge0}$.

\head 1. The totally positive part of $\cb_u$\endhead
\subhead 1.1\endsubhead
Let $G$ be a reductive connected algebraic group over $\CC$ and let $\cb$ be the
variety of Borel subgroups of $G$. 
For $g\in G$ the Springer fibre at $g$ is the subvariety $\cb_g=\{B\in\cb;g\in B\}$ of $\cb$.
The Springer fibres can be quite complicated; they play an important role in representation theory
(for example in character formulas for finite reductive groups over a finite field).
In this paper we assume that a pinning of $G$ is given so that the totally positive part $G_{\ge0}$
of $G$ and the totally positive part $\cb_{\ge0}$ of $\cb$ are defined (see \cite{L94, 2.2, 8.1}). Recall that
$G_{\ge0}$ is a submonoid and a closed subset of $G$ and $\cb_{\ge0}$ is a closed subset of $\cb$
on which $G_{\ge0}$ acts. We are interested in the interaction of the theory of total positivity
with that of Springer fibres. More precisely, for any $g\in G_{\ge0}$ we consider
the closed subset $\cb_{g,\ge0}=\cb_g\cap\cb_{\ge0}$ of $\cb_g$, which we call the totally positive
part of the Springer fibre $\cb_g$. 
Let $\cu$ be the variety of unipotent elements of $G$ and let $\cu_{\ge0}=\cu\cap G_{\ge0}$. 
Our main result is that if $g\in\cu_{\ge0}$, $\cb_{g,\ge0}$ has a surprisingly simple
structure, much simpler than that of $\cb_g$ itself (see Cor. 1.16); namely it has a canonical
cell decomposition which is part of the canonical cell decomposition of $\cb_{\ge0}$. 
By a similar argument we show that an analogous result holds for the fixed point set of $g$ 
on a partial flag manifold intersected with the totally positive part \cite{L98} of that partial 
flag manifold (see 3.9). In \S4 we state a conjecture (see 4.4) on the compatibility of two ways
to define a positive structure on the cells of $\cb_{\ge0}$ and give some of its consequences; this
conjecture extends a result of \cite{L97, \S3} which concerns the open dense cell $\cb_{>0}$ of $\cb_{\ge0}$
defined in \cite{L94, 8.8}. Let $G_{>0}$ the open dense sub-semigroup of $G_{\ge0}$ defined in \cite{L94}.
In \cite{L94, 8.9(c)}, a map $g\m B_g$ from $G_{>0}$ to $\cb_{>0}$
was defined. In \S5 we give a definition of this map which is simpler than that in \cite{L94} and we define
an extension of this map to a map from $G_{\ge0}$ to the set of parabolic subgrups of $G$. As a result we obtain
a new definition of $\cb_{>0}$.

We shall always assume that $G$ is simply laced; the non-simply laced case can be reduced to the
simply laced case by descent.

\subhead 1.2\endsubhead
Recall that we have fixed a pinning of $G$. Thus we are given a maximal torus $T$ of $G$  and a pair 
$B^+,B^-$ of opposed Borel subgroups of $G$ containing $T$ with unipotent radicals $U^+,U^-$. 
The pinning also includes root homomorphisms $x_i:\CC@>>>U^+$,  $y_i:\CC@>>>U^-$ indexed by
a finite set $I$ (corresponding to simple roots). Let $NT$ be the normalizer of $T$ in $G$
and let $W$ be the Weyl group. For $i\in I$ we set $\ds_i=y_i(1)x_i(-1)y_i(1)\in NT$; let $s_i$ be the
image of $\ds_i$ in $W$. Now $W$ is a Coxeter group with simple reflections $\{s_i;i\in I\}$; let $w\m|w|$
be the standard length function. Let $w_I$ be the unique element of maximal length of $W$.
More generally, for any $H\sub I$ we denote by $w_H$ the longest
element in the subgroup $W_H$ of $W$ generated by $\{s_i;i\in H\}$.

Let $\le$ be the standard partial order on $W$.
For $w\in W$ let $\ci_w$ be the set of sequences $(i_1,i_2,\do,i_m)$ in $I$ such that
$w=s_{i_1}s_{i_2}\do s_{i_m}$, $m=|w|$; let $\supp(w)$ be the set of all $i\in I$ which appear in some
(or equivalently any) $(i_1,i_2,\do,i_m)\in\ci_w$. For $w\in W$ we set 
$\dw=\ds_{i_1}\ds_{i_2}\do\ds_{i_m}\in NT$ where $(i_1,i_2,\do,i_m)\in\ci_w$;
this is known to be well defined. 

For $B,B'$ in $\cb$ there is a unique $w\in W$ denoted by $\pos(B,B')$ 
such that for some $g\in G$ we have $gBg\i=B^+,gB'g\i=\dw B^+\dw\i$.
There is a unique isomorphism $\ph:G@>\si>>G$ such that $\ph(x_i(a))=y_i(a)$, $\ph(y_i(a))=x_i(a)$ 
for all $i\in I,a\in\CC$ and $\ph(t)=t\i$ for all $t\in T$. This carries Borel subgroups to Borel subgroups
hence induces an isomorphism $\ph:\cb@>>>\cb$ such that $\ph(B^+)=B^-$, $\ph(B^-)=B^+$. For $i\in I$ we have 
$\ph(\ds_i)=x_i(1)y_i(-1)x_i(1)=\ds_i\i$. Hence $\ph$ induces the identity map on $W$. We show:

(a) {\it Let $B,B'$ in $\cb$ be such that $\pos(B,B')=w$. Then $\pos(\ph(B),\ph(B'))=w_Iww_I$.}
\nl
Assume first that for some $i\in I$ we have $B=B^+,B'=\ds_iB^+\ds_i\i$ so that $\pos(B,B')=s_i$. We have 
$$\align&\pos(\ph(B),\ph(B'))=\pos(B^-,\ds_i\i B^-\ds_i)=\pos(\dw_I B^+\dw_I\i,\ds_i\i\dw_I B^+\dw_I\i\ds_i)\\&
=\pos(B^+,\dw_I\i\ds_i\i\dw_IB^+\dw_I\i\ds_i\dw_I)=w_Is_iw_I.\endalign$$
It follows that (a) holds when $w=s_i$ for some $i\in I$. The general case can be easily reduced to this
special case. This proves (a).

\subhead 1.3\endsubhead
Let $V$ be an irreducible rational $G$-module over $\CC$ 
with a given highest weight vector $\et$.
Let $\b$ be the basis of $V$ (containing $\et$) obtained by 
specializing at $v=1$ the canonical basis \cite{L90} of the corresponding 
module over the corresponding quantized enveloping algebra.
For any $B\in\cb$ let $L_B$ be the unique $B$-stable line in $V$. 
We can assume that $V$ is such that $B\m L_B$ is a bijection

(a) $\cb@>\si>>\cx$
\nl
onto a subset $\cx$ of the set of lines in $V$. Let 

(b) $V_+=\sum_{b\in\b}\RR_{\ge0}b\sub V$.
\nl
If $\x_1\in V_+,\x_2\in V_+$ then $\x_1+\x_2\in V_+$; moreover,

(c) if $\x_1+\x_2=0$, then $\x_1=\x_2=0$.

\subhead 1.4\endsubhead
Following \cite{KL}, for 
any $v,w$ in $W$ we set $\cb_{v,w}=\{B\in B;\pos(B^+,B)=w,\pos(B^-,B)=w_Iv\}$. It is known that
$\cb_{v,w}\ne\emp$ if and only if $v\le w$. The subvarieties $\cb_{v,w}$ with $v\le w$ form a partition
of $\cb$. For $v\le w$ in $W$ we have $ww_I\le vw_I$; moreover,

(a) {\it $\ph$ defines an isomorphism $\cb_{ww_I,vw_I}@>\si>>\cb_{v,w}$.}
\nl
We use that for $B\in\cb$ we have:

$\pos(B^+,B)=vw_I\implies\pos(B^-,\ph(B))=w_Iv$; $\pos(B^-,B)=w_Iww_I\implies\pos(B^+,\ph(B))=w$.

\subhead 1.5\endsubhead
As in \cite{L94} let $G_{\ge0}$ be the submonoid of $G$ generated by $\{x_i(a);i\in I,a\in\RR_{\ge0}\}$,
$\{y_i(a);i\in I,a\in\RR_{\ge0}\}$ and by $\chi(a)$ for various algebraic group
homomorphisms $\chi:\CC^*@>>>T$ and various $a\in\RR_{>0}$; this is
a closed subset of $G$. We have 

(a) $\ph(G_{\ge0})=G_{\ge0}$.
\nl
By \cite{L94, 3.2} we have in the setup of 1.3:

(b) $G_{\ge0}V_+\sub V_+$.
\nl
From \cite{L94, 8.7} it follows that 

(c) $\ph(\cb_{\ge0})=\cb_{\ge0}$.
\nl
Let $\cx_{\ge0}$ be the set of lines in $\cx$ which meet $V_+-\{0\}$. From \cite{L94, 8.17} we have that:

(d) {\it under the bijection 1.3(a), $\cb_{\ge0}$ corresponds to $\cx_{\ge0}$.}
\nl
Using (b),(d) we deduce

(e) $g\in G_{\ge0},B\in\cb_{\ge0}\implies gBg\i\in\cb_{\ge0}$. 
\nl 
(This also follows from \cite{L94, 8.12}.)

\subhead 1.6\endsubhead
Recall that $\cu_{\ge0}=\cu\cap G_{\ge0}$ and that for $u\in\cu_{\ge0}$ we have
$\cb_{u,\ge0}=\cb_u\cap\cb_{\ge0}$. Let $V,V_+$ be as in 1.3.
For $u\in\cu_{\ge0}$ let $V_+^u=\{\x\in V_+;u\x=\x\}$. In particular 
$V_+^{x_i(a)},V_+^{y_i(a)}$ are defined for $i\in I,a\in\RR_{>0}$. We show:

(a) {\it Let $u,u',u''$ in $\cu_{\ge0}$ be such that $u=u'u''$. 
We have $V_+^u=V_+^{u'}\cap V_+^{u''}$.}
\nl
Let $\x\in V_+^{u'}\cap V_+^{u''}$. Then $u\x=u'(u''\x)=u'\x=\x$ hence 
$\x\in V_+^u$. Conversely, let $\x\in V_+^u$. By \cite{L94, 6.2, 3.2} for any $b\in\b$ we have 
$u'b-b\in V_+$, $u''b-b\in V_+$. Hence for any $\x',\x''$ in $V_+$ we have
$u'\x'-\x'\in V_+$, $u''\x''-\x''\in V_+$. Thus $u''\x=\x+\x'$ where $\x'\in V_+$ so that
$\x=u\x=u'(u''\x)=u'(\x+\x')$. We have $u'\x=\x+\x'_1$, $u'\x'=\x'+\x'_2$ where 
$\x'_1\in V_+,\x'_2\in V_+$. Thus $\x=\x+\x'_1+\x'+\x'_2$ so that $\x'_1+\x'+\x'_2=0$; using 1.3(c)
we see that $\x'_1=\x'=\x'_2=0$. 
Thus we have $u''\x=\x,u'\x=\x$ and $\x\in V_+^{u'}\cap V_+^{u''}$. This proves (a).

We show:

(b) {\it Let 
$i\in I,a\in\RR_{>0}$. We have $V_+^{x_i(a)}=V_+^{x_i(1)}$, $V_+^{y_i(a)}=V_+^{y_i(1)}$.}
\nl
Note that $x_i(a)=x_i(1)^a:V@>>>V$ and $x_i(1)=x_i(a)^{1/a}:V@>>>V$. Hence we have 
$\x\in V_+^{x_i(a)}$ if and only 
if $\x\in V_+^{x_i(1)}$. This proves the first equality in (b). The proof of the second equality is entirely 
similar. 

\mpb

Let $u\in\cu_{\ge0}$. From 1.5(d) we see that under the bijection 1.3(a), $\cb_{u,\ge0}$ corresponds to the 
set of lines $L$ in $\cx_{\ge0}$ which are $u$-stable. Since $u:V@>>>V$ is unipotent, for any such 
$L$, $u$ automatically acts on $L$ as identity. Thus we have:

(c) {\it Under 
the bijection 1.3(a), $\cb_{u,\ge0}$ corresponds to the set of lines in $\cx$ which meet $V_+^u-\{0\}$.}

\subhead 1.7\endsubhead
Following \cite{L94, 8.15} we define a partition $\cb_{\ge0}=\sqc_{(v,w)\in W\T W;v\le w}\cb_{\ge0,v,w}$ by
$\cb_{\ge0,v,w}=\cb_{v,w}\cap\cb_{\ge0}$. In \cite{L94, 8.15} it was conjectured that for $v\le w$,

(a) {\it $\cb_{\ge0,v,w}$ is homeomorphic to $\RR_{>0}^{|w|-|v|}$.}
\nl
 This conjecture was proved by Rietsch 
\cite{R98}, \cite{R99} and in a more explicit form by Marsh and Rietsch \cite{MR}.

For $v\le w$ let $[v,w]$ be the subset of $\cx_{\ge0}$ corresponding to $\cb_{\ge0,v,w}$ under the bijection 
$\cb_{\ge0}\lra\cx_{\ge0}$ in 1.5(d). We have $\cx_{\ge0}=\sqc_{(v,w)\in W\T W;v\le w}[v,w]$.

\subhead 1.8\endsubhead
Let $v\le w$ in $W$ and let $\ii=(i_1,i_2,\do,i_m)\in\ci_w$. According to Marsh and Rietsch 
\cite{MR}, there is a unique sequence $t_1,t_2,\do,t_m$  with $t_k\in\{s_{i_k},1\}$
for $k\in[1,m]$, $t_1t_2\do t_m=v$ and such that 
$t_1\le t_1t_2\le\do \le t_1t_2\do t_m$ and $t_1\le t_1s_{i_2},t_1t_2\le t_1t_2s_{i_3},\do,
t_1t_2\do t_{m-1}\le t_1t_2\do t_{m-1}s_{i_m}$.
Following \cite{MR} we define a subset $\cy_{v,w,\ii}$ of $G$ to be the set of products
$g_1g_2\do g_m$ in $G$ where $g_k=\ds_{i_k}$ if $t_k=s_{i_k}$ and
$g_k=y_{i_k}(a_k)$ with $a_k\in\RR_{>0}$ if $t_k=1$; according to \cite{MR}, the map
$\cy_{v,w,\ii}@>>>\cb$, $g\m gB^+g\i$ is a homeomorphism

(a)  $\cy_{v,w,\ii}@>\si>>\cb_{\ge0,v,w}$.
\nl
Moreover we have a homeomorphism

(b) $\cy_{v,w,\ii}@>\si>>\RR_{>0}^{|w|-|v|}$
\nl
(to $g_1g_2\do g_m$ we associate the sequence consisting of the $a_k$ with $k$ such that $t_k=1$.) 
In particular, 
the factors $g_1,g_2,\do,g_m$ are uniquely determined by their product $g_1g_2\do g_m$. The composition
of the inverse of (b) with (a) is a homeomorphism

(c) $\t_\ii:\RR_{>0}^{|w|-|v|}@>\si>>\cb_{\ge0,v,w}$.

\subhead 1.9\endsubhead
We preserve the setup of 1.8. Assume that $s_{i_1}w\le w$. We show:

(a) {\it We have $t_1=1$ if and only if $v\le s_{i_1}w$.}
\nl
If $t_1=1$ then from $v=t_1\do t_m$ we deduce $v=t_2t_3\do t_m$; thus $v$ is equal to a product of
a subsequence of $s_{i_2},s_{i_3},\do, s_{i_m}$ so that $v\le s_{i_1}w$.
Conversely, if $v\le s_{i_1}w$ then by the results in 1.8 applied to $v,s_{i_1}w$ instead of $v,w$,
we can find a sequence  $t'_2,\do,t'_m$ with $t'_k\in\{s_{i_k},1\}$ for $k\in[2,m]$, $t'_2\do t'_m=v$ and 
such that $t'_2\le t'_2t'_3\le\do\le t'_2\do t'_m$ and $t'_2\le t'_2s_{i_3},\do,
t'_2\do t'_{m-1}\le t'_2\do t'_{m-1}s_{i_m}$. Taking $t'_1=1$ we have
$t'_1t'_2\do t'_m=v$,
$t'_1\le t'_1t'_2\le\do\le t'_1t'_2\do t'_m$ and $t'_1\le t'_1s_{i_2},t'_1t'_2\le t'_1t'_2s_{i_3},\do,
t'_1t'_2\do t'_{m-1}\le t'_1t'_2\do t'_{m-1}s_{i_m}$. By uniqueness we have $(t'_1,t'_2,\do,t'_m)=
(t_1,t_2,\do,t_m)$. Thus $t_1=1$ and (a) is proved.

\subhead 1.10\endsubhead
For any $h\in G$ let $[h]:\cb@>>>\cb$ be the map $B\m hBh\i$. 

Let $i\in I,a\in\RR_{>0}$ and let $v,w$ in $W$ be such that $v\le w$. We show:

(a) {\it If 
$v\le s_iw\le w$, then $[y_i(a)]:\cb@>>>\cb$ restricts to a map $\cb_{\ge0,v,w}@>>>\cb_{\ge0,v,w}$ which 
is fixed point free.}

(b) {\it If $s_iw\le w$, $v\not\le s_iw$, then $[y_i(a)]:\cb@>>>\cb$ restricts to the identity map 
$\cb_{\ge0,v,w}@>>>\cb_{\ge0,v,w}$.}

(c) {\it If 
$w\le s_iw$ then $[y_i(a)]:\cb@>>>\cb$ restricts to a map $\cb_{\ge0,v,w}@>>>\cb_{\ge0,v,s_iw}$.}
\nl
Assume first that $s_iw\le w$. We can find $\ii=(i_1,i_2,\do,i_m)\in\ci_w$ such that $i_1=i$. We use the 
notation of 1.8 relative to $v,w,\ii$.

In case (a) we have $t_1=1$, see 1.9(a). In this case, left multiplication by $y_i(a)$ restricts to a map
$\cy_{v,w,\ii}@>>>\cy_{v,w,\ii}$ given by $g_1g_2\do g_m\m g'_1g'_2\do g'_m$ where
$g_1=y_1(a_1),g'_1=y_1(a_1+a)$, $g'_k=g_k$ for $k>1$. To prove (a), it remains to use that $a_1\m a_1+a$ from
$\RR_{>0}$ to $\RR_{>0}$ is fixed point free.

In case (b) we have $t_1=s_i$, see 1.9(a).  In this case, for $g_1g_2\do g_m\in\cy_{v,w,\ii}$ we have
$g_1=\ds_i$ and for some $b\in B^+$ we have
$y_i(a)g_1g_2\do g_m=y_i(a)\ds_i g_2\do g_m=\ds_i x_i(-a)g_2\do g_m=\ds_i g_2\do g_mb$
(the last equality follows by an argument in \cite{MR, 11.9}). Hence (b) holds.

Assume that we are in case (c). Let $B\in\cb_{\ge0,v,w}$, $B'=y_i(a)By_i(a)\i$.
We have $\pos(B^+,B')=\pos(y_i(a)\i B^+y_i(a),B)$. This equals $s_iw$ since
$$\pos(y_i(a)\i B^+y_i(a),B^+)=s_i,\pos(B^+,B)=w$$ and $s_iw\ge w$.
We have $$\pos(B^-,B')=\pos(y_i(a)\i B^-y_i(a),B)=\pos(B^-,B)=v$$ since $y_i(a)\in B^-$.
Thus $B'\in\cb_{v,s_iw}$. By 1.5(e) we have $B'\in\cb_{\ge0}$ hence $B'\in\cb_{\ge0,v,s_iw}$. 
Hence (c) holds.

\subhead 1.11\endsubhead
Let $i\in I,a\in\RR_{>0}$ and let $v,w$ in $W$ be such that $v\le w$. We show:

(a) {\it If 
$v\le s_iv\le w$, then $[x_i(a)]:\cb@>>>\cb$ restricts to a map $\cb_{\ge0,v,w}@>>>\cb_{\ge0,v,w}$ 
which is fixed point free.}

(b) {\it If $v\le s_iv$, $s_iv\not\le w$, then $[x_i(a)]:\cb@>>>\cb$ restricts to the identity map 
$\cb_{\ge0,v,w}@>>>\cb_{\ge0,v,w}$.}

(c) {\it If 
$s_iv\le v$, then $[x_i(a)]:\cb@>>>\cb$ restricts to a map $\cb_{\ge0,v,w}@>>>\cb_{\ge0,s_iv,w}$.}
\nl
We apply 1.10(a)-(c) to $ww_I,vw_I$ instead of $v,w$ (note that $ww_I\le vw_I$) and we apply the automorphism 
$\ph$. We obtain the following statements.

(1) {\it If $ww_I\le s_ivw_I\le vw_I$, then $[x_i(a)]:\cb@>>>\cb$ restricts to a map 
$$\ph(\cb_{\ge0,ww_I,vw_I})@>>>\ph(\cb_{\ge0,ww_I,vw_I})$$ which is fixed point free.}

(2) {\it If $s_ivw_I\le vw_I$, $ww_I\not\le s_ivw_I$ then $[x_i(a)]:\cb@>>>\cb$ restricts to the identity map 
$\ph(\cb_{\ge0,ww_I,vw_I})@>>>\ph(\cb_{\ge0,ww_I,vw_I})$.}

(3) {\it If 
$vw_I\le s_ivw_I$, then $[x_i(a)]:\cb@>>>\cb$ restricts to a map $\ph(\cb_{\ge0,ww_I,vw_I})@>>>
\ph(\cb_{\ge0,ww_I,s_ivw_I})$.}
\nl
It remains to note that 

(d) $\ph(\cb_{\ge0,ww_I,vw_I})=\cb_{\ge0,v,w}$
\nl
for any $v\le w$ in $W$ (see 1.4(a) and 1.5(c)).

\subhead 1.12\endsubhead
Let $i\in I,a\in\RR_{>0}$. From 1.10 we deduce:
$$\{B\in\cb_{\ge0};y_i(a)\in B\}=\sqc_{(v,w)\in W\T W;v\le w,s_iw\le w,v\not\le s_iw}\cb_{\ge0,v,w}.\tag a$$
From 1.11 we deduce:
$$\{B\in\cb_{\ge0};x_i(a)\in B\}=\sqc_{(v,w)\in W\T W;v\le w,v\le s_iv,s_iv\not\le w}\cb_{\ge0,v,w}.\tag b$$

\subhead 1.13\endsubhead
Let $(W\T W)_{disj}=\{(z,z')\in W\T W;\supp(z)\cap\supp(z')=\emptyset\}$. For $(z,z')\in(W\T W)_{disj}$
with $m=|z|,m'=|z'|$, let $\cu_{\ge0,z,z'}$ be the image of the (injective) map
$\RR_{>0}^m\T\RR_{>0}^{m'}@>>>G_{\ge0}$
given by
$$\align&((a_1,a_2,\do,a_m),(a'_1,a'_2,\do,a'_{m'}))\\&\m y_{i_1}(a_1)y_{i_2}(a_2)\do y_{i_m}(a_m)
x_{i'_1}(a'_1)x_{i'_2}(a'_2)\do x_{i'_m}(a'_{m'})\\&
=x_{i'_1}(a'_1)x_{i'_2}(a'_2)\do x_{i'_m}(a'_{m'})y_{i_1}(a_1)y_{i_2}(a_2)\do y_{i_m}(a_m)\endalign$$
where $\ii=(i_1,\do, i_m)\in\ci_w,\ii'=(i'_1,\do, i'_{m'})\in\ci_{w'}$.
This map is injective and its image is independent of the choice of $\ii,\ii'$ (see \cite{L94, 2.7, 2.9}).
Note that $\cu_{\ge0,z,z'}\sub\cu_{\ge0}$. More precisely, we have (see \cite{L94, 6.6}):
$$\cu_{\ge0}=\sqc_{(z,z')\in (W\T W)_{disj}}\cu_{\ge0,z,z'}.$$
We state the following result.

\proclaim{Theorem 1.14} Let $(z,z')\in (W\T W)_{disj}$ and let $u\in\cu_{\ge0,z,z'}$, $J=\supp(z)$,
$J'=\supp(z')$. We have
$$\cb_{u,\ge0}=\cap_{i\in J}\cb_{y_i(1),\ge0}\cap\cap_{j\in J'}\cb_{x_j(1),\ge0}.$$
\endproclaim
Using 1.6(c) we see that the following statement implies the theorem.
$$V_+^u=\cap_{i\in J}V_+^{y_i(1)}\cap \cap_{j\in J'}V_+^{x_j(1)}.\tag a$$
Let $m,m',\ii,\ii'$ be as in 1.13 and let 
$((a_1,a_2,\do,a_m),(a'_1,a'_2,\do,a'_{m'}))$ be corresponding to $u$ as in 1.13.
We argue by induction on $m+m'$. If $m+m'=0$ the result is obvious.
Assume now that $m+m'\ge1$. If $m\ge1$ let $z_1=s_{i_1}z$ and let 
$u_1=y_{i_1}(a_1)\i u$; we have $|z_1|=m-1$ and $u_1\in\cu_{\ge0,z_1,z'},u=y_{i_1}(a_1)u_1$.
If $m'\ge1$ let $z'_1=s_{i'_1}z'$ and let 
$u'_1=x_{i'_1}(a'_1)\i u$; we have $|z'_1|=m'-1$ and $u'_1\in\cu_{\ge0,z,z'_1},u=x_{i'_1}(a'_1)u'_1$.
From 1.6(a),(b), we have
$V_+^u=V_+^{y_{i_1}(1)}\cap V_+^{u_1}$ if $m\ge1$,
$V_+^u=V_+^{x_{i'_1}(1)}\cap V_+^{u'_1}$ if $m'\ge1$.
Since the induction hypothesis is applicable to $V_+^{u_1}$ (if $m\ge1$) and to
$V_+^{u'_1}$ (if $m'\ge1$) we see that (a) is proved. This proves the theorem.

\subhead 1.15\endsubhead
Let $J\sub I,J'\sub I$ be such that $J\cap J'=\emp$. We define 
$$\align&Z_{J,J'}=\{(v,w)\in W\T W;v\le w;s_iw\le w,v\not\le s_iw\qua \frl i\in J;\\&
v\le s_jv,s_jv\not\le w\qua \frl j\in J'\}.\endalign$$
Combining 1.14 with 1.12(a),(b) we obtain our main result. 
\proclaim{Corollary 1.16} Let $(z,z')\in(W\T W)_{disj}$ and let $u\in\cu_{\ge0,z,z'}$,
$J=\supp(z)$, $J'=\supp(z')$. We have
$$\cb_{u,\ge0}=\cup_{(v,w)\in Z_{J,J'}}\cb_{\ge0,v,w}.$$
\endproclaim
Thus $\cb_{u,\ge0}$ admits a canonical, explicit, cell decomposition in which each cell
is part of the canonical cell decomposition of $\cb_{\ge0}$. The zero dimensional cells of
$\cb_{u,\ge0}$ are $\cb_{\ge0,w,w}$ where $w\in W$ is such that the set $\{i\in I;s_iw\le w\}$ contains 
$J$ and is contained in $I-J'$. (For example $w_J$ is such a $w$ since $J\sub I-J'$.) In particular, 
$Z_{J,J'}\ne\emp$, so that $\cb_{u,\ge0}\ne\emp$. (This also follows from \cite{L94, 8.11}.)

\subhead 1.17\endsubhead
In \cite{L94, 2.11}, $G_{\ge0}$ is partitioned into pieces indexed by $W\T W$; by results in \cite{L94}
each piece is a cell. In \cite{L19} the set of pieces of $G_{\ge0}$ is interpreted as a monoid $G(\{1\})$, 
the value of $G$ at the semifield $\{1\}$ 
 with $1$ element, so that pieces appear precisely as the fibres of a natural
surjective map $G_{\ge0}@>>>G(\{1\})=W\T W$ compatible with the monoid structures (the monoid structure on
$W\T W$ thus obtained is not the usual group structure). In particular, the product of two pieces of 
$G_{\ge0}$ is contained in a single piece of $G_{\ge0}$. Similarly, we can view the set of cells
$\cb_{\ge0,v,w}$ of $\cb_{\ge0}$ (see 1.7) as $\cb(\{1\})$ or the value of $\cb$ at the semifield $\{1\}$.
From 1.10 and 1.11 we see that in the action of $G_{\ge0}$ on $\cb_{\ge0}$, the result of applying
a piece of $G_{\ge0}$ to a cell of $\cb_{\ge0}$ is contained in a single cell of $\cb_{\ge0}$. It follows
that the action of $G_{\ge0}$ on $\cb_{\ge0}$ induces an action of $G(\{1\})=W\T W$ on $\cb(\{1\})$. 
We can identify $\cb(\{1\})$ with $\{(v,w)\in W\T W;v\le w\}$. Then the action of $W\T W$ becomes:
$$(s_i,1):(v,w)\m(v,s_i*w), (1,s_i):(v,w)\m(s_i\circ v,w)$$
where for $i\in I, v\in W,w\in W$ we define

$s_i*w=w$ if $s_iw\le w$, $s_i*w=s_iw$ if $w\le s_iw$, 

$s_i\circ v=v$ if $v\le s_iv$, $s_i\circ v=s_iv$ if $s_iv\le v$. 

\mpb

Now let $u\in\cu_{\ge0}$ and let $z,z',J,J'$ be associated to $u$ as in 1.16.
 Let $\cb_u(\{1\})$ be the set of cells of $\cb_{u,\ge0}$ described in 1.16; this set
is in bijection with $Z_{J,J'}$ (see 1.16) and can be viewed as a subset of 
$\cb(\{1\})$.
Let $\fZ_u(\{1\})$ be the set of all $(w,w')\in W\T W=G(\{1\})$ such that
in the $W\T W$-action on $\cb(\{1\})$ (as above), we have $(w,w')\cb_u(\{1\})\sub\cb_u(\{1\})$.
Clearly, $\fZ_u(\{1\})$ is a submonoid of $W\T W=G(\{1\})$.

Let $\fZ(u)$ be the inverse image of $\fZ(u)(\{1\})$ under
the canonical monoid homomorphism $G_{\ge0}@>>>G(\{1\})$ (as above).
Note that $\fZ(u)$ is a submonoid of $G_{\ge0}$  which is
related to the centralizer of $u$ in $G$ (although it is not in general contained in it).
Clearly the $G_{\ge0}$-action of $\cb_{\ge0}$ restricts to a $\fZ(u)$-action on $\cb_{u,\ge0}$.

\subhead 1.18\endsubhead
Let $\ti\cb=\{(u,B)\in\cu\T\cb;u\in B\}$. Let $\ti\cb_{\ge0}=\{(u,B)\in\cu_{\ge0}\T\cb_{\ge0};u\in B\}$.
From 1.16 we can deduce that $\ti\cb_{\ge0}$ is the disjoint union of the sets
$$\ti\cb_{\ge0,z,z',v,w}=\cu_{\ge0,z,z'}\T\cb_{\ge0,v,w}\tag a$$
where $(z,z')$ runs through $(W\T W)_{disj}$ and $(v,w)$ runs through
$Z_{J,J'}$ where $J=\supp(z),J'=\supp(z')$.
 The subset (a) is a cell of dimension $|z|+|z'|+|w|-|v|$. Note that the first projection 
$\ti\cb_{\ge0}@>>>\cu_{\ge0}$ is a trivial fibration over each $\cu_{\ge0,z,z'}$.
Let $\ti\cb(\{1\})$ be the indexing set for the set of cells of $\ti\cb_{\ge0}$. This is the set of
all $(z,z',v,w)\in W^4$ such that $J=\supp(z),J'=\supp(z')$ are disjoint and $(v,w)\in Z_{J,J'}$.

\subhead 1.19\endsubhead
Let $B\in\cb_{\ge0}$ and let $\un B=B\cap\cu_{\ge0}$. This is a closed subset of $B$ closed
under multiplication. (If $u\in\un B,u'\in\un B$ then $u,u'$ are contained in the unipotent radical
of $B$ hence $uu'$ is also contained in that unipotent radical, so that $uu'$ is unipotent.
Since $u\in G_{\ge0},u'\in G_{\ge0}$ we have also $uu'\in G_{\ge0}$. Thus $uu'\in\un B$.) We have
$B\in\cb_{\ge0,v,w}$ for well defined $v\le w$ in $W$. We set
$$\align&\Xi_{v,w}=\{(z,z')\in(W\T W)_{disj};s_iw\le w,v\not\le s_iw\qua \frl i\in\supp(z);\\&
v\le s_jv,s_jv\not\le w\qua \frl j\in\supp(z')\}.\endalign$$
We show:
$$\un B=\sqc_{(z,z')\in\Xi_{v,w}}\cu_{\ge0,z,z'}.\tag a$$
Let $E$ be the right hand side of (a). Let $u\in \un B$. 
Let $(z,z')\in(W\T W)_{disj}$ be such that $u\in\cu_{\ge0,z,z'}$.
Let $J=\supp(z)$, $J'=\supp(z')$. 
We have $B\in\cb_{u,\ge0}$ hence by 1.16 we have $(v,w)\in Z_{J,J'}$ hence $(z,z')\in\Xi_{v,w}$
Thus, $u\in E$. We see that $\un B\sub E$.
Conversely, assume that $u\in E$. We can find $(z,z')\in\Xi_{v,w}$ such that
$u\in \cu_{\ge0,z,z'}$. Let $J=\supp(z)$, $J'=\supp(z')$. We have $(v,w)\in Z_{J,J'}$. Since
$B\in\cb_{\ge0,v,w}$ from 1.16 we see that $B\in\cb_{u,\ge0}$ so that $u\in B$ and $u\in\un B$. We see
that $E\sub\un B$. This proves (a).

\mpb

Let $H=\{i\in I;s_iw\le w,v\not\le s_iw\}$, $H'=\{j\in I;v\le s_jv,s_jv\not\le w\}$.

If $i\in H$ then $v\le w, s_iw\le w$ and this is known to imply $s_iv\le w$. Thus $i\in H\implies i\n H'$
so that 

(b) $H\cap H'=\emp$.
\nl
We have $\Xi_{v,w}=\{(z,z')\in(W\T W)_{disj};\supp(z)\sub H,\supp(z')\sub H'\}$. Using (b) we see that 

(c) $\Xi_{v,w}=\{(z,z')\in W\T W;\supp(z)\sub H,\supp(z')\sub H'\}$.
\nl
Hence (a) becomes:
$$\un B=\sqc_{(z,z')\in W\T W;\supp(z)\sub H,\supp(z')\sub H'}\cu_{\ge0,z,z'}.\tag d$$
From (d) we see that $\un B$ has a canonical cell decomposition with cells indexed by (c). One of these cells
is $\cu_{\ge0,w_H,w_{H'}}$ (of dimension $|w_H|+|w_{H'}|$); all other cells have dimension strictly less than
$|w_H|+|w_{H'}|$ and are contained in the closure of $\cu_{\ge0,w_H,w_{H'}}$. We see that 

(e) {\it $\un B$ is connected of dimension $|w_H|+|w_{H'}|$.}
\nl
If $(v,w)=(1,w_I)$ then $H=H'=\emp$ and from (d) we see that $\un B=\cu_{\ge0,1,1}=\{1\}$.
Recall that $\cb_{\ge0,1,w_I}=\cb_{>0}$. We see that:

(f) {\it If $u\in\cu_{\ge0}$ and $B\in\cb_u$ satisfies $B\in\cb_{>0}$ then $u=1$.}
\nl
It is likely that, more generally, for any $B\in\cb_{>0}$, $B\cap G_{\ge0}$ consists of semisimple
elements.

\subhead 1.20\endsubhead
Let $g\in G_{\ge0}$ and let $[g]:\cb@>>>\cb$ be as in 1.10. Let $\cb_{g,\ge0}$ be as in 1.1. 
If $v\le w$ in $W$ then, by 1.17, there are three possibilities:

(i) $[g]\cb_{\ge0,v,w}=\cb_{\ge0,v',w'}$ where $v'\le w',(v',w')\ne(v,w)$;

(ii) $[g]\cb_{\ge0,v,w}=\cb_{\ge0,v,w}$ and $B\m[g]B$, $\cb_{\ge0,v,w}@>>>\cb_{\ge0,v,w}$ has empty fixed point set;

(iii) $[g]\cb_{\ge0,v,w}=\cb_{\ge0,v,w}$ and $B\m[g]B$, $\cb_{\ge0,v,w}@>>>\cb_{\ge0,v,w}$ has non-empty fixed point set (denoted by
$\cb_{g,\ge0,v,w}$).
\nl
It follows that $\cb_{g,\ge0}=\sqc_{v,w}\cb_{g,\ge0,v,w}$ where $v,w$ is as in (iii).
We conjecture that if $v,w$ is as in (iii), then $\cb_{g,\ge0,v,w}$ is a cell. (When $g$ is unipotent this holds by 1.16.)

\head 2. Examples\endhead
\subhead 2.1\endsubhead
Assume that $G=SL_{n+1}(\CC)$ and $I=\{1,2,\do,n\}$ with $n\ge2$ and with 
$s_1s_2,s_2s_3,\do,s_{n-1}s_n$ of order $3$. Let $J=\{1,2,\do,n-1\}$. Let $-$ be the unit element of $W$. 
We have
$(w_J,-)\in(W\T W)_{disj}$ and $Z_{J,\emp}$ consists of   
$$\align&(w_J,w_J),(w_Js_n,w_Js_n),(w_J,w_Js_n),(w_Js_ns_{n-1},w_Js_ns_{n-1}),(w_Js_n,w_Js_ns_{n-1}),\\&
\do,(w_Js_ns_{n-1}\do s_2s_1,w_Js_ns_{n-1}\do s_2s_1),(w_Js_ns_{n-1}\do s_2,w_Js_ns_{n-1}\do s_2s_1).\tag a
\endalign$$
Thus if $u\in\cu_{\ge0,w_J,-}$ then $\cb_{u,\ge0}$ is a union of $n+1$ cells of dimension $0$ and $n$ 
cells of dimension $1$. 

We have also $(w_J,s_n)\in(W\T W)_{disj}$ and $Z_{J,\{n\}}$ consists of the pairs in (a) other than
the last two. 
Thus, if $u\in\cu_{\ge0,z,s_n}$ then $\cb_{u,\ge0}$ is a union of $n$ cells of dimension $0$ and $n-1$ cells of 
dimension $1$. 

Assume now that $n=3$. We have $(s_1s_3,-)\in (W\T W)_{disj}$. Let $J=\{1,3\}$. Then $Z_{J,\emp}$ consists of
$$\align&(132,132),(13,132),(1321,1321),(132,1321),(1323,1323),(132,1323),\\&
(13231,13231),\\&(1321,13213),(1323,13213),(132,13213),(132312,132312),(13231,132132).\tag b\endalign$$
(We write $i_1i_2\do$ instead of $s_{i_1}s_{i_2}\do$.)
Thus, if $u\in\cu_{\ge0,s_1s_3,-}$ then $\cb_{u,\ge0}$ is a union of $6$ cells of dimension $0$, $6$ cells of 
dimension $1$ and one cell of dimension $2$.

We have also $(s_1s_3,s_2)\in (W\T W)_{disj}$ and $Z_{J,\{2\}}$ consists of the pairs in (b) other than the 
last two. Thus, if $u\in\cu_{\ge0,s_1s_3,s_2}$ then $\cb_{u,\ge0}$ is a union of $5$ cells of dimension $0$, 
$5$ cells of dimension $1$ and one cell of dimension $2$.

In these examples $\cb_{u,\ge0}$ is contractible; but in the last two
examples, $\cb_{u,\ge0}$ is not of pure dimension, unlike $\cb_u$.

We have $(s_1,-)\in (W\T W)_{disj}$. Then $Z_{\{1\},\emp}$ consists of
$$\align&(2,2),(21,21),(2,21),(23,23),(2,23),(212,212),(21,212),\\&
(232,232),(23,232),(213,213),(21,213),(23,213),(2,213),(2123,2123),\\&
(212,2123),(213,2123),(21,2123),(2321,2321),(232,2321),(231,2321),\\&
(23,2321),(2132,2132),(213,2132),(232,2132),(212,2132),(21,2132),(23,2132),\\&
(21232,21232),(2123,21232),(2132,21232),(212,21232),(21,21232),\\&
(23212,23212),(2321,23212),(2312,23212),(232,23212),(231,23212),\\&
(23,23212),(213213,213213),(32132,213213),(12312,213213),(2132,213213).\endalign$$
(We write $i_1i_2\do$ instead of $s_{i_1}s_{i_2}\do$.)
Thus if $u\in\cu_{\ge0,s_1,-}$ then
$\cb_{u,\ge0}$ is a union of cells of dimension $\le3$, two of which have dimension $3$.

In each of the examples above, for any cell $\cb_{\ge0,v,w}$ of maximal dimension of $\cb_{u,\ge0}$, we have
$w=vw_H$ where $H\sub I$ is such that $u$ is conjugate in $G$ to a regular unipotent element in the 
subgroup of $G$ generated by $\{y_i(a),x_i(a);i\in I-H,a\in\CC\}$ and by $T$. This is likely to be a 
general phenomenon.

For such $u$ one can show that $|w_H|\le\dim\cb_u$ (complex dimension); this is compatible
with $|w_H|=\dim\cb_{u,\ge0}$ (real dimension).

\subhead 2.2\endsubhead
In this subsection we assume that $G,I$ are as in 2.1 and $n=2$.
In this case we can take $V$ in 1.3 to be the adjoint representation of $G$.
The canonical basis $\b$ of $V$ can be denoted by $X_{-12},X_{-1},X_{-2},t_1,t_2,X_1,X_2,X_{12}$
and the action of $x_i(a),y_i(a),i\in I,a\in\CC$ is as follows:

$x_i(a)X_{12}=X_{12}$

$x_i(a)X_j=X_j+aX_{12}$ if $i\ne j\in I$

$x_i(a)X_j=X_j$ if $i=j$

$x_i(a)X_{-j}=X_{-j}$ if $i\ne j\in I$

$x_i(a)X_{-j}=X_{-j}+at_j+a^2X_j$ if $i=j$

$x_i(a)X_{-12}=X_{-12}+aX_{-j}$ if $i\ne j\in I$

$x_i(a)t_j=t_j+aX_i$ if $i\ne j\in I$

$x_i(a)t_j=t_j+2aX_i$ if $i=j$

$y_i(a)X_{12}=X_{12}+aX_j$ if $i\ne j\in I$

$y_i(a)X_j=X_j$ if $i\ne j\in I$

$y_i(a)X_j=X_j+at_j+a^2X_{-j}$ if $i=j$

$y_i(a)X_{-j}=X_{-j}+aX_{-12}$ if $i\ne j\in I$

$y_i(a)X_{-j}=X_{-j}$ if $i=j$

$y_i(a)X_{-12}=X_{-12}$

$y_i(a)t_j=t_j+aX_{-i}$ if $i\ne j\in I$

$y_i(a)t_j=t_j+2aX_{-i}$ if $i=j$.

The set $\cx_{\ge0}$ consists of all
$$a_{-12}X_{-12}+a_{-1}X_{-1}+a_{-2}X_{-2}+c_1t_1+c_2t_2+a_1X_1+a_2X_2+a_{12}X_{12}$$
where 
$a_{-12},a_{-1},a_{-2},c_1,c_2,a_1,a_2,a_{12}$ are in $\bold R_{\ge0}$ (not all $0$) such that 
$$\align&a_2a_{-12}=c_2a_{-1},a_1a_{-12}=c_1a_{-2},a_{-1}a_{12}=c_1a_2,\\&
a_{-2}a_{12}=c_2a_1,a_{12}(c_1+c_2)=a_1a_2,a_{-12}(c_1+c_2)=a_{-1}a_{-2},\\&
c_1c_2=a_{12}a_{-12},c_1(c_1+c_2)=a_1a_{-1},c_2(c_1+c_2)=a_2a_{-2}\tag a\endalign$$
modulo the homothety action of $\bold R_{>0}$.

The subsets $[v,w]$ of $\cx_{\ge0}$ can be described
as follows (the coefficients $$a_{-12},a_{-1},a_{-2},c_1,c_2,a_1,a_2,a_{12}$$ are required to be in 
$\bold R_{>0}$ and are taken up to simultaneous multiplication by an element in $\bold R_{>0}$):

$[121,121]$:     $\{a_{-12}X_{-12}\},$

$[12,12]$:  $\{a_{-1}X_{-1}\},$

$[21,21]$:  $\{a_{-2}X_{-2}\},$

$[2,2]$: $\{a_1X_1\},$

$[1,1]$: $\{a_2X_2\},$

$[-,-]$: $\{a_{12}X_{12}\},$

$[21,121]$: $\{a_{-12}X_{-12}+a_{-2}X_{-2}\},$

$[12,121]$: $\{a_{-12}X_{-12}+a_{-1}X_{-1}\},$

$[1,12]$: $\{a_{-1}X_{-1}+a_2X_2\},$

$[2,21]$: $\{a_{-2}X_{-2}+a_1X_1\},$

$[-,2]$: $\{ a_1X_1+a_{12}X_{12}\},$

$[-,1]$: $\{a_2X_2+a_{12}X_{12}\},$

$[2,12]$: $\{a_{-1}X_{-1}+c_1t_1+a_1X_1; a_{-1}a_1=c_1^2\},$

$[1,21]$: $\{a_{-2}X_{-2}+c_2t_2+a_2X_2; a_{-2}a_2=c_2^2\},$
$$\align&[2,121]:\qua \{a_{-12}X_{-12}+a_{-1}X_{-1}+a_{-2}X_{-2}+c_1t_1+a_1X_1 , a_{-1}a_1=c_1^2,\\&
c_1a_{-12}=a_{-1}a_{-2}, a_1a_{-12}=c_1a_{-2}\},\endalign$$
$$\align& [1,121]:\qua  
\{a_{-12}X_{-12}+a_{-1}X_{-1}+a_{-2}X_{-2}+c_2t_2+a_2X_2;a_{-2}a_2=c_2^2,\\&c_2a_{-12}=a_{-1}a_{-2}, 
a_2a_{-12}=c_2a_{-1}\},\endalign$$
$$\align&[-,12]:\qua
\{a_{-1}X_{-1}+c_1t_1+a_1X_1+a_2X_2+a_{12}X_{12};a_{-1}a_1=c_1^2,\\&c_1a_{12}=a_1a_2,a_{-1}a_{12}=c_1a_2\},
\endalign$$
$$\align:&[-,21]: \{a_{-2}X_{-2}+c_2t_2+a_2X_2+a_1X_1+a_{12}X_{12};a_{-2}a_2=c_2^2,\\&
c_2a_{12}=a_1a_2, a_{-2}a_{12}=c_2a_1\},\endalign$$ 
$$\align&[-,121]:\\& \{a_{-12}X_{-12}+a_{-1}X_{-1}+a_{-2}X_{-2}+c_1t_1+c_2t_2+a_1X_1+a_2X_2+a_{12}X_{12}\\&
\text{ such that (a) holds}\}.\endalign$$ 

\subhead 2.3\endsubhead
We preserve the setup of 2.2.

(a) {\it For any $v\le w$ in $W$ there is a well defined subset $[[v,w]]$ of $\b$ such that $[v,w]$ consists 
of all lines in $\cx_{\ge0}$ which contain some vector spanned by an $\RR_{>0}$-linear combination of 
vectors in $[[v,w]]$.}
\nl
We have

$[[121,121]]=\{X_{-12}\},$

$[[12,12]]=\{X_{-1}\},$

$[[21,21]]=\{X_{-2}\},$

$[[2,2]]=\{X_1\},$

$[[1,1]]=\{X_2\},$

$[[-,-]]=\{X_{12}\},$

$[[21,121]]=\{X_{-12},X_{-2}\},$

$[[12,121]]=\{X_{-12},X_{-1}\},$

$[[1,12]]=\{X_{-1},X_2\},$

$[[2,21]]=\{X_{-2},X_1\},$

$[[-,2]]=\{X_1,X_{12}\},$

$[[-,1]]=\{X_2,X_{12}\},$

$[[2,12]]=\{X_{-1},t_1,X_1\},$

$[[1,21]]=\{X_{-2},t_2,X_2\},$

$[[2,121]]=\{X_{-12},X_{-1},X_{-2},t_1,X_1\},$

$[[1,121]=\{X_{-12},X_{-1},X_{-2},t_2,X_2\},$

$[[-,12]]=\{X_{-1},t_1,X_1,X_2,X_{12}\},$

$[[-,21]]=\{X_{-2},t_2,X_2,X_1,X_{12}\},$
 
$[[-,121]]=\{X_{-12},X_{-1},X_{-2},t_1,t_2,X_1,X_2,X_{12}\}$.

We have 

(b) {\it  $[[v,w]]=[[v,121]]\cap[[-,w]]$ for any $v\le w$ in $W$.}
\nl
(c) {\it There is a well defined partition $\b=\sqc_{z\in W}\b_z^-$ such that
$$[[-,w]]=\sqc_{z\in W;z\le w}\b_z^-$$
for any $w\in W$. There is a well defined partition $\b=\sqc_{z\in W}\b_z^+$ such that}
$$[[v,121]]=\sqc_{z\in W;v\le z}\b_z^+.$$
We have 

$\b^-_{-}=\{X_{12}\}$, $\b^-_1=\{X_2\}$, $\b^-_2=\{X_1\}$, $\b^-_{12}=\{X_{-1},t_1\}$,
$\b^-_{21}=\{X_{-2},t_2\}$, $\b^-_{121}=\{X_{-12}\}$.

$\b^+_{-}=\{X_{12}\}$, $\b^+_1=\{t_2,X_2\}$, $\b^+_2=\{t_1,X_1\}$, 
$\b^+_{21}=\{X_{-2}\}$, $\b^+_{12}=\{X_{-1}\}$, $\b^+_{121}=\{X_{-12}\}$.

\subhead 2.4\endsubhead
We return to the general case. 
Now 2.3(a),(b),(c) make sense in the general case (in (b),(c) we replace $121$ by $w_I$ and $-$ by the
unit element of $W$); we expect that these statements hold in the general case. 
In particular $\b^-_z\sub \b$ and $\b^+_z\sub\b$ are defined for $z\in W$. 

Let $\et^-$ be the unique vector in $\b$ such that the
stabilizer of $\CC\et^-$ in $G$ is $B^-$. We can regard $V$ naturally as a module over the
universal enveloping algebra of the Lie algebra of $G$ hence as a module
over the universal enveloping algebra $\fU^+$ of the Lie algebra of $U^+$
and as a module over the universal enveloping algebra $\fU^-$ of the Lie algebra of $U^-$.
Now $\fU^+$ (resp. $\fU^-$) has a canonical basis $\ti\b^+$ (resp. $\ti\b^-$), see \cite{L90}) and 
the map $c^+:\ti b\m\ti b\et^-$ (resp. $c^-:\ti b\m\ti b\et$) from $\ti\b^+$ (resp. $\ti\b^-$) to $V$
has image $\b\cup\{0\}$. 
From \cite{L19, 10.2} we have a partition $\ti\b^+=\sqc_{w\in W}\ti\b^+_w$; similarly we have a partition
$\ti\b^-=\sqc_{w\in W}\ti\b^-_w$. We expect that for $z\in W$, $\b^-_z$ is equal to $c^-(\ti\b^-_z)$ with 
$0$ removed and that $\b^+_z$ is equal to $c^+(\ti\b^+_{zw_I})$ with $0$ removed. 

\head 3. Partial flag manifolds\endhead
\subhead 3.1\endsubhead
We fix $H\sub I$. Let $W^H$ be the set of all $w\in W$ such that $w$ has minimal length in $wW_H$.
Let $P_H$ be the subgroup of $G$ generated by 
$\{x_i(a);i\in I,a\in\CC\}$, $\{y_i(a);i\in H,a\in\CC\}$ and by $T$ (a parabolic subgroup
containing $B^+$). Let $\cp_H$ be the variety whose points are the subgroups of $G$ conjugate to
$P_H$. (We have $P_\emp=B^+,\cp_\emp=\cb$.) Define  $\p_H:\cb@>>>\cp_H$ by $B\m P$ where $P\in\cp_H$
contains $B$. As observed in \cite{L98} to any $P\in\cp_H$ we can attach two Borel subgroups $B',B''$ of 
$P$ such that $\pos(B^+,B')=b\in W^H,\pos(B^-,B'')=a'\in W^H$; moreover $B',B''$ (hence $a',b$) are 
uniquely determined by $P$. Let $c=\pos(B',B'')\in W_H$ and let $a=w_Ia'$. We have $a\le bc$ 
(since $B''\in\cb_{a,bc}$) and $ac\i\le b$ (since $B'\in\cb_{ac\i,b}$). Conversely, if 

(a) $(a,b,c)\in (w_IW^H)\T W^H\T W_H$ satisfy $a\le bc$ or equivalently $ac\i\le b$,
\nl
then the set of all $P\in\cp_H$ which give rise as above to $a,b,c$ is non-empty; we denote this set by
 $\cp_{H,a,b,c}$. The subsets $\cp_{H,a,b,c}$ form a partition of $\cp_H$ indexed by the set 
$\cp_H(\{1\})$ of triples $(a,b,c)$ as in (a).

Following \cite{L98} we set $\cp_{H,\ge0}=\p_H(\cb_{\ge0})$, $\cp_{H,>0}=\p_H(\cb_{>0})$. 
For $(a,b,c)\in\cp_H(\{1\})$ we set $$\cp_{H,\ge0,a,b,c}:=\cp_{H,a,b,c}\cap\cp_{H,\ge0}$$ so that we have
$$\cp_{H,\ge0}=\sqc_{(a,b,c)\in\cp_H(\{1\})}\cp_{H,\ge0,a,b,c}.$$
Let 

$\cz=\{(r,t)\in(w_IW^H)\T W;r\le t\}$, $\cz'=\{(r',t')\in W\T W^H;r'\le t'\}$. 
\nl
We have bijections 

$\a:\cp_H(\{1\})@>\si>>\cz$, $(a,b,c)\m(a,bc)$, $\a':\cp_H(\{1\})@>\si>>\cz'$, $(a,b,c)\m(ac\i,b)$.
\nl
We shall write ${}_{r,t}\cp_{H,\ge0}$ instead of $\cp_{H,\ge0,a,b,c}$ where $(r,t)\in\cz,(a,b,c)=\a\i(r,t)$
 and ${}^{r',t'}\cp_{H,\ge0}$ instead of $\cp_{H,\ge0,a,b,c}$ where $(r',t')\in\cz',(a,b,c)=\a'{}\i(r',t')$.
Thus we have 

$\cp_{H,\ge0}=\sqc_{(r,t)\in\cz}{}_{r,t}\cp_{H,\ge0}=\sqc_{(r',t')\in\cz'}{}^{r',t'}\cp_{H,\ge0}$
\nl
and 

${}_{r,t}\cp_{H,\ge0}={}^{r',t'}\cp_{H,\ge0}$ if $(r',t')=\a'\a\i(r,t)$.
\nl
In \cite{R98, p.50,51} it is shown that for $(r,t)\in\cz$, $\p_H$ restricts to a bijection 

(b) $\ti\a:\cb_{\ge0,r,t}@>\si>>{}_{r,t}\cp_{H,\ge0}$
\nl
and that for $(r',t')\in\cz'$, $\p_H$ restricts to a bijection 

(c) $\ti\a':\cb_{\ge0,r',t'}@>\si>>{}^{r',t'}\cp_{H,\ge0}$.
\nl
This implies (by 1.7(a)) that $\cp_{H,\ge0,a,b,c}$ is a cell for any 
$(a,b,c)\in\cp_H(\{1\})$, so that the various $\cp_{H,\ge0,a,b,c}$ form a 
cell decomposition of $\cp_{H,\ge0}$; this justifies the notation $\cp_H(\{1\})$. 

From 1.5(d) we deduce:

(d) $g\in G_{\ge0},P\in\cp_{H,\ge0}\implies gPg\i\in\cp_{H,\ge0}$. 
\nl
Let $\cp=\sqc_{H\sub I}\cp_H$, $\cp_{\ge0}=\sqc_{H\sub I}\cp_{H,\ge0}$, $\cp_{>0}=\sqc_{H\sub I}\cp_{H,>0}$.

\subhead 3.2\endsubhead
Let $H$ be as in 3.1. Let $V_H$ be an irreducible rational $G$-module over $\CC$ such that $\{g\in G;g\cl=\cl\}=P_H$
for some (necessarily unique) line $\cl$ in $V$; let $\et\in\cl-\{0\}$. 
Let $\b$ be the basis of $V_H$ (containing $\et$) obtained by 
specializing at $v=1$ the canonical basis \cite{L90} of the corresponding 
module over the corresponding quantized enveloping algebra.
Let $\cx_H$ be the $G$-orbit of $\cl$ in the set of lines in $V_H$. We have a bijection 

(a) $\cp_H@>\si>>\cx_H$
\nl
given by $P\m L_P$ where $L_P$ is the unique $P$-stable line in $V_H$. Let 

(b) $V_{H+}=\sum_{b\in\b}\RR_{\ge0}b\sub V_H$.
\nl
From \cite{L98} we see that $V_H$ above can be chosen so that:

(c) {\it under the bijection (a), $\cp_{H,\ge0}$ corresponds to the set of lines in $\cx_H$ which meet $V_{H+}-\{0\}$.}
\nl
In the sequel we assume that $V_H$ has been chosen so that (c) holds.

\subhead 3.3\endsubhead
For $u\in\cu_{\ge0}$ let $\cp_{H,u,\ge0}=\{P\in\cp_{H,\ge0};u\in P\}$, $V_{H+}^u=\{\x\in V_{H+};u\x=\x\}$. 
The proof of (a),(b),(c) below is identical to that of 1.6(a),(b),(c).

(a) {\it Let $u,u',u''$ in $\cu_{\ge0}$ be such that $u=u'u''$. We have $V_{H+}^u=V_{H+}^{u'}\cap V_{H+}^{u''}$.}

(b) {\it Let $i\in I,a\in\RR_{>0}$. We have $V_{H+}^{x_i(a)}=V_{H+}^{x_i(1)}$, $V_{H+}^{y_i(a)}=V_{H+}^{y_i(1)}$.}

(c) {\it Under the bijection 3.2(c), $\cp_{H,u,\ge0}$ corresponds to the set of lines in $\cx_H$ which meet $V_{H+}^u-\{0\}$.}

\subhead 3.4\endsubhead
For any $h\in G$ let $[h]_H:\cp_H@>>>\cp_H$ be the map $P\m hPh\i$. 

Let $i\in I,a\in\RR_{>0}$. Let $(r,t)\in\cz$. Now (a),(b),(c) below follow immediately from 1.10(a),(b),(c) using 3.1(b).

(a) {\it If $r\le s_it\le t$, then $[y_i(a)]_H:\cp_H@>>>\cp_H$ restricts to a map ${}_{r,t}\cp_{H,\ge0}@>>>{}_{r,t}\cp_{H,\ge0}$ which 
is fixed point free.}

(b) {\it If $s_it\le t$, $r\not\le s_it$, then $[y_i(a)]_H:\cp_H@>>>\cp_H$ restricts to the identity map ${}_{r,t}\cp_{H,\ge0}@>>>{}_{r,t}\cp_{H,\ge0}$.}

(c) {\it If $t\le s_it$ then $[y_i(a)]_H:\cp_H@>>>\cp_H$ restricts to a map 
${}_{r,t}\cp_{H,\ge0}@>>>{}_{r,s_it}\cp_{H,\ge0}$; note that $(r,s_it)\in\cz$.}

\subhead 3.5\endsubhead
Let $i\in I,a\in\RR_{>0}$. Let $(r',t')\in\cz'$. Now (a),(b),(c) below follow immediately from 1.11(a),(b),(c) using 3.1(c).

(a) {\it If $r'\le s_ir'\le t'$, then $[x_i(a)]_H:\cp_H@>>>\cp_H$ restricts to a map ${}^{r',t'}\cp_{H,\ge0}@>>>{}^{r',t'}\cp_{H,\ge0}$ which 
is fixed point free.}

(b) {\it If $r'\le s_ir'$, $s_ir'\not\le t'$, then $[x_i(a)]_H:\cp_H@>>>\cp_H$ restricts to the identity map 
${}_{r',t'}\cp_{H,\ge0}@>>>{}_{r',t'}\cp_{H,\ge0}$.}

(c) {\it If $s_ir'\le r'$ then $[x_i(a)]_H:\cp_H@>>>\cp_H$ restricts to a map ${}^{r',t'}\cp_{H,\ge0}@>>>{}^{s_ir',t'}\cp_{H,\ge0}$; note that 
$(s_ir',t')\in\cz'$.}

\subhead 3.6\endsubhead
Let $i\in I,a\in\RR_{>0}$. From 3.4 we deduce:
$$\{P\in\cp_{H,\ge0};y_i(a)\in P\}=\sqc_{(r,t)\in\cz;s_it\le t,r\not\le s_it}({}_{r,t}\cp_{H,\ge0}).\tag a$$
From 3.5 we deduce:
$$\{P\in\cp_{H,\ge0};x_i(a)\in P\}=\sqc_{(r',t')\in\cz';r'\le s_ir',s_ir'\not\le t}{}^{r',t'}\cp_{H,\ge0}.\tag b$$

The proof of the following result is entirely similar to that of 1.14, using 3.3(a),(b) instead of 1.6(a),(b).
\proclaim{Theorem 3.7} Let $(z,z')\in (W\T W)_{disj}$ and let $u\in\cu_{\ge0,z,z'}$, $J=\supp(z)$,
$J'=\supp(z')$. We have
$$\cp_{H,u,\ge0}=\cap_{i\in J}\cp_{H,y_i(1),\ge0}\cap\cap_{j\in J'}\cp_{H,x_j(1),\ge0}.$$
\endproclaim

\subhead 3.8\endsubhead
Let $J\sub I,J'\sub I$ be such that $J\cap J'=\emp$. We set
$$\align&Z_{H;J,J'}=\{((r,t),(r',t'))\in\cz\T\cz';\\& 
(r',t')=\a'\a\i(r,t),s_it\le t,r\not\le s_it\qua \frl i\in J;
r'\le s_jr',s_jr'\not\le t' \qua\frl j\in J'\}.\endalign$$
Combining 3.7 with 3.6(a),(b) we obtain the following result. 
\proclaim{Corollary 3.9} Let $(z,z')\in(W\T W)_{disj}$ and let $u\in\cu_{\ge0,z,z'}$,
$J=\supp(z)$, $J'=\supp(z')$. We have
$$\cp_{H,u,\ge0}=\cup_{((r,t),(r',t'))\in Z_{H;J,J'}}({}_{r,t}\cp_{H,\ge0}).$$
\endproclaim
Thus $\cp_{H,u,\ge0}$ admits a canonical, explicit, cell decomposition in which each cell
is part of the canonical cell decomposition of $\cp_{H,\ge0}$.
The zero dimensional cells of $\cp_{H,u,\ge0}$ are indexed by
$$\{((r,r),(r',r'));r'\in W^H,r=r'w_H, s_ir\le r\qua \frl i\in J,r'\le s_jr'\qua \frl j\in J'\}.$$
The last set contains for example $((r,r),(r',r'))$ where $$r'=w_Jw_{J\cap H},r=w_Jw_{J\cap H}w_H.$$

\subhead 3.10\endsubhead
From 3.4 and 3.5 we see that in the action of $G_{\ge0}$ on $\cp_{H,\ge0}$, the result of applying
a piece of $G_{\ge0}$ to a cell of $\cp_{H,\ge0}$ is contained in a single cell of $\cp_{H,\ge0}$. It follows
that the action of $G_{\ge0}$ on $\cp_{H,\ge0}$ induces an action of $G(\{1\})=W\T W$ on $\cp_H(\{1\})$. 
We can identify $\cp_H(\{1\})$ with 
$$\{((r,t),(r',t'))\in\cz\T\cz'; (r',t')=\a'\a\i(r,t)\}.$$
Then the action of $W\T W$ becomes:
$$(s_i,1):((r,t),(r',t'))\m((r,s_i*t),\a'\a\i(r,s_i*t)),$$     
$$(1,s_i):((r,t),(r',t'))\m(\a\a'{}\i(s_i\circ r',t'),(s_i\circ r',t')),$$    
(notation of 1.17).

\subhead 3.11\endsubhead
Let $P\in\cp_{H,\ge0}$ and let $\un P=P\cap\cu_{\ge0}$. This is a closed subset of $P$.
We have $P\in {}_{r,t}\cp_{H,\ge0}={}^{r',t'}\cp_{H,\ge0}$ for a well defined
$((r,t),(r',t'))\in\cz\T\cz'$ such that $(r',t')=\a'\a\i(r,t)$.
We set
$$\align&\Xi_{H,(r,t),(r',t')}=\{(z,z')\in(W\T W)_{disj};s_it\le t,r\not\le s_it\qua \frl i\in\supp(z);\\&
r'\le s_jr',s_jr'\not\le t' \qua\frl j\in \supp(z')\}.\endalign$$
We show:
$$\un P=\sqc_{(z,z')\in\Xi_{H,(r,t),(r',t')}}\cu_{\ge0,z,z'}.\tag a$$
Let $E$ be the right hand side of (a). Let $u\in\un P$. 
Let $(z,z')\in(W\T W)_{disj}$ be such that $u\in\cu_{\ge0,z,z'}$.
Let $J=\supp(z)$, $J'=\supp(z')$. 
We have $P\in\cp_{H,u,\ge0}$ hence by 3.9 we have $((r,t),(r',t'))\in Z_{H;J,J'}$
hence $(z,z')\in\Xi_{H,(r,t),(r',t')}$. Thus, $u\in E$. We see that $\un P\sub E$.
Conversely, assume that $u\in E$. We can find $(z,z')\in\Xi_{H,(r,t),(r',t')}$ such that
$u\in \cu_{\ge0,z,z'}$. Let $J=\supp(z)$, $J'=\supp(z')$. 
We have $((r,t),(r',t'))\in Z_{H;J,J'}$. Since $P\in{}_{r,t}\cp_{H,\ge0}$, 
from 3.9 we see that $P\in\cp_{H,u,\ge0}$ so that $u\in P$ and $u\in\un P$. We see
that $E\sub\un P$. This proves (a).

From (a) we see that $\un P$ has a canonical cell decomposition with cells indexed by $\Xi_{H,(r,t),(r',t')}$.

\head 4. A conjecture and its consequences\endhead
\subhead 4.1\endsubhead
An $\RR_{>0}$-positive structure on a set $\fX$ is a finite collection $f_e:\RR_{>0}^m@>\si>>\fX$ of 
bijections ($e\in E$) with $m\ge0$ fixed, such that for any $e,e'$ in $E$, 
$f_e\i f_{e'}:\RR_{>0}^m@>>>\RR_{>0}^m$ is 
admissible in the sense of \cite{L19, 1.2}. If $(\fX;f_e:\RR_{>0}^m@>>>\fX, e\in E)$,
$(\ti\fX;\ti f_{\ti e}:\RR_{>0}^{\ti m}@>>>\ti\fX,\ti e\in\ti E)$ 
are two sets with $\RR_{>0}$-positive structure, a map $\x:\fX@>>>\ti\fX$ is said to be a morphism if
for some (or equivalently any) $e\in E,\ti e\in\ti E$, the map 
$\ti f_{\ti e}\i\x f_e:\RR_{>0}^m@>>>\RR_{>0}^{\ti m}$ is admissible in the sense of \cite{L19, 1.2}. 
We say that $\x$ is an isomorphism if it is a bijective morphism and $\x\i$ is a morphism.

If $\fX,\fX'$ 
have positive $\RR_{>0}$-structures, then $\fX\T\fX'$ has a natural $\RR_{>0}$-positive structure.

\subhead 4.2\endsubhead
Let $T_{>0}=T\cap G_{\ge0}$. For any basis $\c_*=(\c_1,\c_2,\do,\c_n)$ of 
$\Hom(\CC^*,T)$ we have a bijection $\RR_{>0}^n@>\si>>T_{>0}$ given by 
$(a_1,\do,a_n)\m\c_1(a_1)\do\c_n(a_n)$. These bijections (for various $\c_*$) define
an $\RR_{>0}$-positive structure on $T_{>0}$.

Let $v\le w$ in $W$. From \cite{R08} we see that the bijections 
$\t_\ii:\RR_{>0}^{|w|-|v|}@>\si>>\cb_{\ge0,v,w}$ (see 1.8(c))
with $\ii\in\ci_w$ form an $\RR_{>0}$-positive structure 
on $\cb_{\ge0,v,w}$. 

From 1.10 we see that for $i\in I$, 

(a) {\it the map $\RR_{>0}\T\cb_{\ge0,v,w}@>>>\cb_{\ge0,v,s_i*w}$ given by 
$(a,B)\m y_i(a)By_i(a)\i$ is a well defined morphism of sets with $\RR_{>0}$-positive structure.}

From the definitions, for any $\chi\in\Hom(\CC^*,T)$, 

(b) {\it the map $\RR_{>0}\T\cb_{\ge0,v,w}@>>>\cb_{\ge0,v,w}$ given by 
$(a,B)\m \chi(a)B\chi(a)\i$ is a well defined morphism of sets with $\RR_{>0}$-positive structure.}
 
\subhead 4.3\endsubhead
Replacing in 4.2 $v,w$ by $w_Iw,w_Iv$, we deduce that the bijections
$$\t_{\ii'}:\RR_{>0}^{|w|-|v|}=\RR_{>0}^{|w_Iv|-|w_Iw|}@>\si>>\cb_{\ge0,w_Iw,w_Iv}$$
 with $\ii'\in\ci_{w_Iv}$ form an $\RR_{>0}$-positive structure on $\cb_{\ge0,w_Iw,w_Iv}$ and that
for $i\in I$, the map 
$$\RR_{>0}\T\cb_{\ge0,w_Iw,w_Iv}@>>>\cb_{\ge0,w_Iw,s_i*(w_Iv)}$$ 
given by $(a,B)\m y_i(a)By_i(a)\i$ is a well defined morphism of sets with $\RR_{>0}$-positive structure. 
We state:

\proclaim{Conjecture 4.4} The bijection $\cb_{\ge0,ww_I,vw_I}@>\si>>\cb_{\ge0,v,w}$ defined by $\ph$
(see 1.11(d)) is an isomorphism of $\RR_{>0}$-positive structures.
\endproclaim
When $v=1,w=w_I$ this can be deduced from \cite{L97, \S3}.

{\it In the remainder of this section we assume that this conjecture holds.}

\subhead 4.5\endsubhead
Let $v,w$ be as in 4.2. Using 4.4 we can reformulate the last statement in 4.3 as follows. For $i\in I$, 

(a) {\it the map $\RR_{>0}\T\cb_{\ge0,v,w}@>>>\cb_{\ge0,s_i\circ v,w}$ given by 
$(a,B)\m x_i(a)Bx_i(a)\i$ is a well defined morphism of sets with $\RR_{>0}$-positive structure. }

\subhead 4.6\endsubhead
We consider a sequence of admissible maps (see \cite{L19, 1.2}):
$$\align&\Ph_1:\RR_{>0}\T\RR_{>0}^{n_0}@>>>\RR_{>0}^{n_1}, \Ph_2:\RR_{>0}\T\RR_{>0}^{n_1}@>>>\RR_{>0}^{n_2},
\\&\do,\Ph_\s:\RR_{>0}\T\RR_{>0}^{n_{\s-1}}@>>>\RR_{>0}^{n_\s}\endalign$$
where $n_0,n_1,\do,n_\s$ are in $\NN$. We define $\Ph:\RR_{>0}^\s\T\RR_{>0}^{n_0}@>>>\RR_{>0}^{n_\s}$ by
$$\Ph((a_1,a_2,\do,a_\s),b)=\do\Ph_3(a_3,\Ph_2(a_2,\Ph_1(a_1,b)))\do)$$
where $b\in \RR_{>0}^{n_0}$. Clearly, $\Ph$ is admissible.

\subhead 4.7\endsubhead
Consider the inverse image $G_{r,-s}$ of $(r,s)\in W\T W$ under the map $G_{\ge0}@>>>W\T W$ in 1.17.
Now $G_{r,-s}$ has a natural $\RR_{>0}$-positive structure (see \cite{L19}). 
Let $(i_1,i_2,\do,i_m)\in\ci_r,(j_1,j_2,\do,j_l)\in\ci_s$.
Let $v\le w$ be elements of $W$. Define $v'\le w'$ by 
$v'=s_{i_1}\circ(\do s_{i_{m-1}}\circ(s_{i_m}\circ v))\do)$,
$w'=s_{j_1}*(\do s_{j_{l-1}}*(s_{j_l}\circ w))\do)$ (notation of 1.17).
The $G_{\ge0}$-action on $\cb_{\ge0}$ restricts to a map $G_{r,-s}\T\cb_{\ge0,v,w}@>>>\cb_{\ge0,v',w'}$.

(a) {\it This is a morphism of sets with $\RR_{>0}$-positive structure.}
\nl
Indeed, our map can be identified with a $\Ph$ in 4.6 where each of $\Ph_1,\do,\Ph_\s$ 
is a map as in 4.2(a),(b) or 4.5(a) hence is admissible.

\subhead 4.8\endsubhead
We now fix a semifield $K$. Let $(\fX;f_e:\RR_{>0}^m@>>>\fX, e\in E)$ be a set with an $\RR_{>0}$-positive 
structure. For any $e,e'$ in $E$ the admissible bijection $f_e\i f_{e'}:\RR_{>0}^m@>\si>>\RR_{>0}^m$ induces 
a bijection $\un f_{e,e'}:K^m@>\si>>K^m$. (This is obtained by replacing the indeterminates which appear in
the formula for $f_e\i f_{e'}$ by elements of $K$ instead of elements of $\RR_{>0}$.) There is a well 
defined set $\fX(K)$ with
bijections $\un f_e:K^m@>\si>>\fX(K)$ such that $\un f_{e,e'}=\un f_e\i\un f_{e'}$ for any $e,e'$ in $E$.
If $\x:\fX@>>>\fX'$ is a morphism of sets with $\RR_{>0}$-positive structure then $\x$ induces a map
of sets $\x(K):\fX(K)@>>>\fX'(K)$.

\subhead 4.9\endsubhead
In the setup of 4.7, the sets $G_{r,-s}(K),\cb_{\ge0,v,w}(K),\cb_{\ge0,v',w'}(K)$ are defined as in
4.8 and by 4.8(a), the map $G_{r,-s}\T\cb_{\ge0,v,w}@>>>\cb_{\ge0,v',w'}$ induces a map

(a) $G_{r,-s}(K)\T\cb_{\ge0,v,w}(K)@>>>\cb_{\ge0,v',w'}(K)$.
\nl
We set $G(K)=\sqc_{(r,s)\in W\T W}G_{r,-s}(K)$, $\cb(K)=\sqc_{(v,w)\in W\T W;v\le w}\cb_{\ge0,v,w}(K)$.
The maps (a) define a map $G(K)\T \cb(K)@>>>\cb(K)$. This is an action of $G(K)$ (with the monoid
structure induced from that of $G_{\ge0}$) on the set $\cb(K)$.

\subhead 4.10\endsubhead
Let $H\sub I$. For any $(a,b,c)\in \cp_H(\{1\})$, the cell $\cp_{H,\ge0,a,b,c}$ in $\cp_{H,\ge0}$ has
a $\RR_{>0}$-positive structure via a  bijection as in 3.1(a) or as in 3.1(b) (these two bijections
define the same $\RR_{>0}$-positive structure, as we see easily from the definitions). Thus
the set $\cp_{H,\ge0,a,b,c}(K)$ is defined. We set 
$$\cp_H(K)=\sqc_{(a,b,c)\in\cp_H(\{1\})}\cp_{H,\ge0,a,b,c}(K).$$
An argument similar to that in 4.9 shows that the $G_{\ge0}$-action on $\cp_{H,\ge0}$ (see 3.1(d)) induces 
an action of the monoid $G(K)$ on $\cp_H(K)$.

\subhead 4.11\endsubhead
Let $u\in\cu_{\ge0}$. Let $z,z',J,J'$ be as in 1.16. We set 
$$\cb_u(K)=\cup_{(v,w)\in Z_{J,J'}}\cb_{\ge0,v,w}(K),$$
$$\fZ(u)(K)=\cup_{(r,s)\in \fZ_u(\{1\})}G_{r,-s}(K).$$
(notation of 1.15, 1.17). Then the action of $G(K)$ on $\cb(K)$ restricts to an action 
of $\fZ(u)(K)$ (with the monoid structure induced from that of $G(K)$) on the set $\cb_u(K)$.

\head 5. The map $g\m P_g$ from $G_{\ge0}$ to $\cp_{\ge0}$\endhead
\subhead 5.1\endsubhead
In this section we give a new, simpler, definition of the map $g\m B_g$ from $G_{>0}$ to $\cb_{>0}$
in \cite{L94, 8.9(c)} and we extend it to a map $g\m P_g$ from $G_{\ge0}$ to $\cp_{\ge0}$.

\subhead 5.2\endsubhead
For a closed subgroup $G'$ of $G$ we denote by $\fL G'$ the Lie algebra of $G'$. Let $\fg=\fL G$. 
Let $g\in G_{\ge0}$. We associate to $g$ a parabolic subgroup $P=P_g$ of $G$ containing $g$ as follows.
By \cite{L19, 9.1(a)}, we have $\fg=\op_{a\in\RR_{>0}}\fg_a$ where $\fg_a$ is the generalized $a$-eigenspace of $\Ad(g):\fg@>>>\fg$.
It follows that we have $\fg=\fg_{<1}\op\fg_1\op\fg_{>1}$ where $\fg_{<1}=\op_{a:0<a<1}\fg_a$, $\fg_{>1}=\op_{a:a>1}\fg_a$ and
$\fg_{\le1}=\fg_{<1}\op\fg_1,\fg_{\ge1}=\fg_1\op\fg_{>1}$ are opposed parabolic subalgebras of $\fg$ with common Levi subalgebra $\fg_1$.
Let $P=P_g$ be the parabolic subgroup of $G$ with $\fL P=\fg_{\ge1}$.
Let $L=L_g$ be the Levi subgroup of $P$ with $\fL L=\fg_1$. Note that $L$ is the centralizer in $G$ of the semisimple part $g_s$ of $g$.
In particular we have $g\in L$; more precisely, $g$ is a central element of $L$ times a unipotent element of $L$.
For example, if $g\in G_{\ge0}$ is unipotent then $P_g=G$.

\subhead 5.3\endsubhead
We now assume that $g\in G_{>0}$. We show:

(a) $B_g=P_g$.
\nl
Let $U^+_{>0}\sub U^+,U^-_{>0}\sub U^-$ be as in \cite{L94, 2.12}. Let $T_{>0}\sub T$ be as in 4.7.
By \cite{L94, 8.10} we can find $u\in U^+_{>0},u'\in U^-_{>0},t\in T_{>0}$ such that 
$g=u'tuu'{}\i$ and all eigenvalues of $\Ad(t)$ on $\fL U^+$ are $>1$. Then all eigenvalues of $\Ad(t)$ on $\fL B^+$ are $\ge1$.
Now $tu$ is $U^+$-conjugate to $t$ hence all eigenvalues of $\Ad(tu)$ on $\fL B^+$ are $\ge1$. It follows that all eigenvalues of 
$\Ad(g)$ on $\fL(u'B^+u'{}\i)$ are $\ge1$, so that $\fL(u'B^+u'{}\i)\sub\fL P_g$. Since $t$ (and $tu$) is regular semisimple,
$L_g$ is a maximal torus of $G$ so that $P_g\in\cb$; this implies that $\fL(u'B^+u'{}\i)=\fL P_g$. From the definition of $\cb_{>0}$
we have $u'B^+u'{}\i\in=\cb_{>0}$. Since $g\in u'B^+u'{}\i$, we have $u'B^+u'{}\i=B_g$. This proves (a). 

\subhead 5.4\endsubhead
Let $g\in G_{\ge0}$. We show:

(a) $P_g\in\cp_{\ge0}$.
\nl
By \cite{L94, 4.4} we can find a sequence $g_1,g_2,g_3,\do$ in $G_{>0}$ such that 
$g=\lim_{n\to\iy}g_n$ in $G$. Since $\cb$ is compact, some subsequence of the sequence of Borel subgroups $P_{g_n}$ 
converges in $\cb$ to a Borel subgroup $B\in\cb$. We can assume that $\lim_{n\to\iy}P_{g_n}=B$ in $\cb$. Since 
$P_{g_n}\in\cb_{>0}\sub\cb_{\ge0}$ and $\cb_{\ge0}$ is closed in $\cb$, we have $B\in\cb_{\ge0}$.
Let $\fp_n=\fL P_{g_n}$, $n=1,2,\do$, $\fp=\fL P_g$, and let $\fb=\fL B$. We show:

(b) $\fb\sub\fp$.
\nl
We have $\lim_{n\to\iy}\fp_n=\fb$ in the Grassmannian of $\dim\fb$-dimensional subspaces of $\fg$.
Since $\fp_n$ is stable under $\Ad(g_n)$ for $n=1,2,\do$, we see that $\fb$ is stable under $\Ad(g)$.
Moreover since all eigenvalues of $\Ad(g_n):\fp_n@>>>\fp_n$ are $\ge1$, we see that all
eigenvalues of $\Ad(g):\fb@>>>\fb$ are $\ge1$. Hence (b) holds.

From (b) we deduce that $B\sub P_g$. Since $B\in\cb_{\ge0}$ we see that $P_g\in\cp_{\ge0}$; this proves (a) and completes the
verification of the statements in 5.1.

\subhead 5.5\endsubhead
Let $B\in\cb_{>0}$. We show:

(a) $B\cap G_{>0}$ is non-empty.
\nl
We have $B=uB^+u\i$ for some $u\in U^-_{>0}$. Let $\ti u\in U^+_{>0}$. By \cite{L94, 7.2} we can find 
$t\in T_{>0}$ such that $u\ti u tu\i\in G_{>0}$. We have $u\ti u tu\i\in uB^+u\i=B$. Thus
$u\ti u tu\i\in B\cap G_{>0}$ and (a) follows.

\mpb

We conjecture that $B\cap G_{>0}$ is homeomorphic to $\RR_{>0}^{\dim B}$.
Consider for example the case where $G=SL_2(\CC)$. In this case there is a unique $z\in\RR_{>0}$ such that
$B$ equals $\{\left(\sm a& b\\c&d\esm\right);a,b,c,d\text{ in }\CC,ad=bc+1,c+dz=z(a+bz)\}$.
Then $B\cap G_{>0}$ can be identified with the set $\{(a,b,c,d)\in\RR_{>0}^4;ad=bc+1,c+dz=z(a+bz)\}$ or (via the sustitution
$d=(bc+1)/a$), with the set $X_z=\{(a,b,c)\in\RR_{>0}^3;c+z\fra{bc+1}{a}=z(a+bz)\}=\{(a,b,c)\in\RR_{>0}^3;(az-c)(bz+a)=z\}$.
Setting $\e=a-cz\i$ we can identify $X_z$ with 
$X'_z=\{(b,c,\e)\in\RR_{>0}^3;\e^2+\e(bz+cz\i)-1=0\}$. The map $X'_z@>>>\RR_{>0}^2, (b,c,\e)\m(b,c)$ is a homeomorphism:
for any $(b,c)\in\RR_{>0}^2$ there is a unique $\e\in\RR_{>0}$ such that $\e^2+\e(bz+cz\i)-1=0$. 
This proves the conjecture in our case.

\subhead 5.6\endsubhead
We show:

(a) {\it If $g\in G_{>0}$ then $P_g\in\cb$; the image of the map $g\m P_g$ from $G_{>0}$ to $\cb$ is exactly $\cb_{>0}$.}
\nl
The fact that $P_g\in\cb_{>0}$ for $g\in G_{>0}$ follows from 5.3(a). Let $B\in\cb_{>0}$. Let $g\in B\cap G_{>0}$ (see 5.5(a)). 
From \cite{L94, 8.9(a)} we see that $B=B_g$ hence by 5.3(a) we have $B=P_g$. This proves (a).

\mpb

Now (a) provides a new definition of $\cb_{>0}$; it is the image of the map $g\m P_g$ from $G_{>0}$ to $\cb$.

Let $G_{reg}$ be the set of regular elements in $G$ and let $G_{reg,\ge0}=G_{reg}\cap G_{\ge0}$.
We conjecture that the map $g\m P_g$ from $G_{\ge0}$ to $\cp_{\ge0}$ is surjective and that, moreover, its 
restriction $G_{reg,\ge0}@>>>\cp_{\ge0}$ is surjective.

\enddocument